\documentclass[a4paper,12pt]{amsart}
\usepackage{amssymb}
\usepackage[francais]{babel}
\input xypic
\xyoption{all}
\xyoption{arc}

\def\theorem#1{{\refstepcounter{theo}\label{#1}\noindent\bf Th\'eor\`eme 
\arabic{section}.\arabic{theo}.}}
\def\proposition#1{{\refstepcounter{theo}\label{#1}\noindent\bf Proposition 
\arabic{section}.\arabic{theo}.}}
\def\lemma#1{{\refstepcounter{theo}\label{#1}\noindent\bf Lemme  
\arabic{section}.\arabic{theo}.}}
\def\scholie#1{{\refstepcounter{theo}\label{#1}\noindent\bf   
(\arabic{section}.\arabic{theo}) \quad}}

\def\remark#1{{\refstepcounter{theo}\label{#1}\noindent\sc Remarque
\arabic{section}.\arabic{theo} - }}

\def\equat{\refstepcounter{theo}$$~}
\def\endequat{\leqno{\boldsymbol{(\arabic{section}.\arabic{theo})}}~$$}

\newcounter{soussection}[section]
\def\soussection#1{\refstepcounter{soussection}
\noindent{\bf \arabic{section}.\Alph{soussection}. #1.}}

\newcounter{numero}[section]

\def\CM{{\mathbb{C}}}

\def\QM{{\mathbb{Q}}}

\def\ZM{{\mathbb{Z}}}

\def\AG{{\mathfrak A}}

\def\SG{{\mathfrak S}}

\def\a{\alpha}
\def\b{\beta}

\def\g{\gamma}

\def\e{\varepsilon}
\def\ph{\varphi}

\def\L{\Lambda}
\def\m{\mu}

\def\r{\rho}
\def\s{\sigma}

\def\t{\tau}

\def\z{\zeta}

\def\BC{{\mathcal{B}}}

\def\WC{{\mathcal{W}}}

\def\Gb{{\mathbf G}}

\def\Lb{{\mathbf L}}

\def\Ob{{\mathbf O}}

\def\Sb{{\mathbf S}}
\def\Tb{{\mathbf T}}

\def\ib{{\mathbf i}}

\def\nb{{\mathbf n}}

\def\pb{{\mathbf p}}

\def\Mti{{\tilde{M}}}

\def\Sti{{\tilde{S}}}

\def\pht{{\tilde{\varphi}}}

\def\mut{{\tilde{\mu}}}

\def\sigt{{\tilde{\s}}}

\def\diag{\mathop{\mathrm{diag}}\nolimits}

\def\Id{\mathop{\mathrm{Id}}\nolimits}

\def\Ker{\mathop{\mathrm{Ker}}\nolimits}

\def\Tr{\mathop{\mathrm{Tr}}\nolimits}

\def\tete#1{\par\leavevmode\makebox[0.7cm]{$(\mathrm{#1})$}}

\def\to{\rightarrow}
\def\longto{\longrightarrow}

\def\longmapright#1{\hspace{0.3em}\smash{
     \mathop{\longrightarrow}\limits^{#1}}\hspace{0.3em}}

\def\fonction#1#2#3#4#5{\begin{array}{rccc}
{#1} : & {#2} & \longto & {#3} \\
& {#4} & \longmapsto & {#5} 
\end{array}}

\def\incl{\hspace{0.05cm}{\subset}\hspace{0.05cm}}

\def\DS{\displaystyle}
\def\SS{\scriptstyle}

\def\fin{~$\SS \blacksquare$}
\def\finl{~$\SS \square$}

\def\matrice#1{\left(\begin{array}{ccccccccccccccccccc}#1\end{array}\right)}

\def\lexp#1#2{\kern\scriptspace\vphantom{#2}^{#1}\kern-\scriptspace#2}
\def\le{\hspace{0.1em}\mathop{\leqslant}\nolimits\hspace{0.1em}}
\def\ge{\hspace{0.1em}\mathop{\geqslant}\nolimits\hspace{0.1em}}

\mathchardef\lllllll="3278
\def\SEC{$\lllllll$}
\mathchardef\inferieur="321E
\mathchardef\superieur="321F
\def\arobas{\char'100}

\def\eqna{\begin{eqnarray*}}
\def\endeqna{\end{eqnarray*}}

\def\cad{c'est-\`a-dire }

\def\irr{irr\'eductible }

\def\resp{respectivement }

\def\www{{\bigwedge}^2}
\def\wgras{\b_\wedge}
\def\ssss{{\sqrt{\Sb\Lb}}}
\def\ISO{{\ib\Sb\Ob}}

\begin{document}

\begin{centerline}{\Large \bf Une (nouvelle ?) construction du groupe}\end{centerline}

\medskip

\begin{centerline}{\Large \bf de r\'eflexion complexe 
${\boldsymbol{G_{31}}}$}\end{centerline}

\bigskip

\begin{centerline}{\sc C. Bonnaf\'e\footnote{CNRS - UMR 6623, Laboratoire de 
Math\'ematiques de Besan\c{c}on, 16 Route de Gray, 25030 BESAN\c{C}ON Cedex, FRANCE, 
{\tt bonnafe\arobas math.univ-fcomte.fr}}}\end{centerline}

\medskip

\begin{centerline}{31 octobre 2002}\end{centerline}

\bigskip

\bigskip

\begin{quotation}
{\small \noindent {\bf R\'esum\'e :} Dans la classification de 
Shephard-Todd des groupes de r\'eflexions, le groupe not\'e $G_{31}$ 
appara\^{\i}t comme le seul groupe de r\'eflexion irr\'eductible 
en rang $4$ qui ne peut pas \^etre engendr\'e par $4$ r\'eflexions. 
Nous en donnons ici une construction totalement \'el\'ementaire 
\`a partir du groupe de Weyl de type $B_6$.}
\end{quotation}

\bigskip

\bigskip

Si $Z$ d\'esigne le sous-groupe central d'ordre $2$ du groupe 
de r\'eflexion complexe $G_{31}$ (voir la classification de 
Shephard-Todd \cite{ST}), alors $G_{31}/Z \not\simeq W(D_6)$ (o\`u $W(D_6)$ 
d\'esigne un groupe de Weyl de type $D_6$), bien que l'on ait 
une suite exacte
$$1 \longto (\ZM/2\ZM)^5 \longto G_{31}/Z \longto \SG_6 \longto 1.$$
En fait, $G_{31}/Z(G_{31}) \not\simeq W(D_6)/Z(W(D_6))$ 
(voir \cite[table 3]{BMR}).

Cependant, on peut v\'erifier gr\^ace \`a {\tt GAP} que 
$D(G_{31})/Z \simeq D(W(D_6))$. En voici une explication~: 
$G_{31}$ agit sur un espace vectoriel $V$ de dimension $4$ 
donc $G_{31}/Z$ agit sur la puissance ext\'erieure deuxi\`eme 
$\www V$ qui est de dimension $6$. C'est \`a travers cette 
action que l'on obtient l'isomorphisme annonc\'e. 

En fait, \`a partir de ce constat, 
nous proposons ici de renverser le point de vue. 
Nous nous donnons un groupe de Weyl $W_6$ de type $B_6$ 
sur $\www V$ (notons que $W_6$ contient comme sous-groupe distingu\'e 
d'indice $2$ un groupe de Weyl de type $D_6$) et construisons 
\`a partir de lui et de fa\c{c}on \'el\'ementaire 
un sous-groupe de $\Gb\Lb(V)$ dont nous montrons que ce ne 
peut-\^etre que le groupe $G_{31}$. Nous d\'eduisons 
des propri\'et\'es classiques de $W_6$ et de ses sous-groupes 
(par exemple de l'existence d'un automorphisme non int\'erieur 
de $\SG_6$) certaines propri\'et\'es du groupe $G_{31}$ 
(structure du groupe d\'eriv\'e, de son plus grand $2$-sous-groupe 
distingu\'e, nombre minimal de r\'eflexions 
n\'ecessaires pour engendrer $G_{31}$...).

\medskip

Dans la preuve de nos r\'esultats, nous n'utilisons pas {\tt GAP}. 
C'est tout de m\^eme un peu artificiel~: il nous a \'et\'e d'un 
grand secours  pour nous donner des id\'ees de d\'emonstration. 
Il nous a aussi servi \`a \'etablir certaines des remarques 
que nous avons rajout\'ees au long de ce texte.

\bigskip

\noindent{\sc Notations - } Si $H$ est un groupe fini, nous notons 
$D(H)$ son groupe d\'eriv\'e, $Z(H)$ son centre et $O_2(H)$ son 
plus grand $2$-sous-groupe distingu\'e. Si $g$ et $h$ sont deux 
\'el\'ements de $H$, nous posons $[g,h]=ghg^{-1}h^{-1}$. 
Si $n$ est un entier 
naturel non nul, nous notons $\SG_n$ (\resp $\AG_n$) le 
groupe sym\'etrique (\resp altern\'e) de degr\'e $n$. 

\bigskip

\section{Construction de $G_{31}$}~

\medskip

\soussection{Puissance ext\'erieure deuxi\`eme d'un espace 
de dimension $4$} Soit $V$ un espace vectoriel complexe de 
dimension $4$. Nous fixons une fois pour toutes un g\'en\'erateur 
$\e$ de ${\bigwedge}^4 V$. Le choix de ce g\'en\'erateur nous permet 
d'identifier ${\bigwedge}^4 V$ avec $\CM$ et donc 
de construire une forme bilin\'eaire
$$\fonction{\wgras}{\www V \times \www V}{\CM}{(x,y)}{x \wedge y.}$$
Il est imm\'ediat que $\wgras$ est sym\'etrique et non d\'eg\'en\'er\'ee.

Nous notons $\L : \Gb\Lb(V) \to \Gb\Lb(\www V)$, $g \mapsto g \wedge g$. 
C'est un morphisme de groupes alg\'ebriques. On a
\equat\label{ker}
\Ker \L=\{\Id_V , - \Id_V\}.
\endequat 
Le noyau de $\L$ sera not\'e $Z$ dans la suite. 
D'autre part, si $g \in \Gb\Lb(V)$ et si $x$, $y \in \www V$, alors
\equat\label{det}
\det \L(g) = (\det g)^3
\endequat
et
\equat\label{sym}
\wgras(\L(g)(x),\L(g)(y)) = (\det g) \wgras(x,y).
\endequat
La derni\`ere \'egalit\'e est simplement une autre \'ecriture 
de la d\'efinition du d\'eterminant. Les \'egalit\'es \ref{det} et 
\ref{sym} montrent que l'image de $\Sb\Lb(V)$ est contenue 
dans $\Sb\Ob(\www V,\wgras)$. 
En fait, pour des raisons de dimension et de connexit\'e, on a 
$$\L(\Sb\Lb(V))=\Sb\Ob(\www V,\wgras).$$
Par cons\'equent, $\L$ induit un isomorphisme de groupes alg\'ebriques 
\equat\label{iso}
\Sb\Lb(V)/\{\Id_V,-\Id_V\} \simeq \Sb\Ob(\www V, \wgras).
\endequat

\bigskip

\remark{a3d3} L'isomorphisme \ref{iso} explique l'\'egalit\'e des diagrammes 
de Dynkin de type $A_3$ et $D_3$. Il montre aussi que 
$\Sb\Lb(V) \simeq \Sb\pb\ib\nb(\www V, \wgras)$, \cad que 
$\Sb\Lb_4(\CM) \simeq \Sb\pb\ib\nb_6(\CM)$.\finl

\bigskip

\remark{b2c2} Soit $\psi^*$ une forme bilin\'eaire altern\'ee 
non d\'eg\'en\'er\'ee sur $V$. On peut voir $\psi^*$ comme un \'el\'ement 
de $\www V^*$. Or, $\www V^* \simeq (\www V)^*$ et 
$\wgras$ induit un isomorphisme $(\www V)^* \simeq \www V$. 
Notons $\psi$ l'\'el\'ement de $\www V$ correspondant \`a $\psi^*$ 
via ces isomorphismes. Alors $\Sb\pb(V,\psi^*)$ stabilise $\psi$, 
donc stabilise $\psi^\perp$ (l'orthogonal pour la forme bilin\'eaire 
sym\'etrique $\wgras$). 
On obtient donc un morphisme de groupes alg\'ebriques 
$\L' : \Sb\pb(V,\psi^*) \to \Sb\Ob(\psi^\perp,\wgras)$. On a 
$\Ker \L'=\Ker \L$ et $\L'$ est surjectif pour 
des raisons de dimension et de connexit\'e. On a donc 
construit un isomorphisme de groupes alg\'ebriques
$$\Sb\pb(V,\psi^*)/\{\Id_V,-\Id_V\} \simeq \Sb\Ob(\psi^\perp,\wgras).$$
Cet isomorphisme explique l'\'egalit\'e des diagrammes 
de Dynkin de type $B_2$ et $C_2$ et montre que 
$\Sb\pb(V,\psi^*) \simeq \Sb\pb\ib\nb(\psi^\perp,\wgras)$, \cad que 
$\Sb\pb_4(\CM) \simeq \Sb\pb\ib\nb_5(\CM)$.\finl

\bigskip

\soussection{Groupes de Weyl de type $B_6$ et $D_6$ sur $\www V$} 
Nous fixons maintenant une base orthonormale $\BC=(\e_1,\e_2,\e_3,\e_4,\e_5,\e_6)$ 
de $\www V$. Nous notons $\Tb$ le groupe des automorphismes de l'espace 
vectoriel $\www V$ dont la matrice dans la base $\BC$ est diagonale. 
Quand nous \'ecrirons un \'el\'ement de $\Gb\Lb(\www V)$ sous forme 
matricielle, il sera sous-entendu que c'est relativement \`a la base $\BC$. 
Nous identifions le groupe $\SG_6$ avec le sous-groupe de $\Ob(\www V, \wgras)$ 
des permutations des \'el\'ements de la base $\BC$. Notons 
que $\SG_6$ normalise $\Tb$ et que l'application 
$$\fonction{\g}{\Tb \rtimes \SG_6}{\CM^\times}{t \rtimes \s}{\det t}$$
est un morphisme de groupes alg\'ebriques. On note 
$$\pi : \Tb \rtimes \SG_6 \to \SG_6$$ 
la projection canonique. Posons maintenant 
\eqna
W_6&=&\Ob(\www V, \wgras) \cap (\Tb \rtimes \SG_6),\\
W_6^+&=&W_6 \cap \Ker(\det),\\
W_6^\prime &=& W_6 \cap \Ker \g,\\
\WC_6&=&<W_6^+ , i \Id_{\www V} > \\
{\mathrm{et}}\hskip4.9cm\WC_6^\prime &=& \WC_6 \cap \Ker \g.\hskip7cm
\endeqna
Ici, $i$ d\'esigne un nombre complexe tel que $i^2=-1$. 
Alors $W_6$ (\resp $W_6^\prime$) est un groupe de Weyl de type 
$B_6$ (\resp $D_6$). On pose maintenant 
$$A_6=W_6 \cap \Tb=\{\diag(a_1,a_2,a_3,a_4,a_5,a_6) \in \Tb~|~\forall 1 \le k \le 6,~
a_k^2=1\}$$ 
$$A_6^\prime = W_6^\prime \cap A_6 = \{\diag(a_1,a_2,a_3,a_4,a_5,a_6) \in A_6~|~
a_1a_2a_3a_4a_5a_6=1\}.\leqno{\mathrm{et}}$$
Alors $A_6 \simeq (\ZM/2\ZM)^6$ et $A_6^\prime \simeq (\ZM/2\ZM)^5$. 
On v\'erifie facilement le r\'esultat suivant~:

\bigskip

\scholie{eng} {\it $A_6^\prime$ est engendr\'e par 
$([\s,a])_{\s \in \AG_6,~a \in A_6^\prime}$.}

\bigskip

Il d\'ecoule de \ref{eng} que $D(D(W_6))=D(W_6) = A_6^\prime \rtimes \AG_6$. 
D'autre part, on a 
$$A_6^\prime=W_6^+ \cap A_6=\WC_6 \cap A_6 
= \WC_6^\prime \cap A_6=D(W_6) \cap A_6$$
et on a des suites exactes
\equat\label{w}
1 \longto A_6 \longto W_6 \longmapright{\DS{\pi}} \SG_6 \longto 1,
\endequat
\equat\label{w'}
1 \longto A_6^\prime \longto W_6^\prime \longmapright{\DS{\pi}} \SG_6 \longto 1,
\endequat
\equat\label{w+}
1 \longto A_6^\prime \longto W_6^+ \longmapright{\DS{\pi}} \SG_6 \longto 1,
\endequat
\equat\label{wc}
1 \longto \WC_6 \cap \Tb \longto \WC_6 \longmapright{\DS{\pi}} \SG_6 \longto 1,
\endequat
\equat\label{wc'}
1 \longto A_6^\prime \longto \WC_6^\prime \longmapright{\DS{\pi}} \SG_6 \longto 1
\endequat
et
\equat\label{dw}
1 \longto A_6^\prime \longto D(W_6) \longmapright{\DS{\pi}} \AG_6 \longto 1.
\endequat
Par cons\'equent, on a  
$$D(W_6)=D(W_6^\prime)=D(W_6^+)=D(\WC_6)=D(\WC_6^\prime)=W_6^\prime \cap W_6^+,$$
$$|W_6|=|\WC_6|=46080=2^6. 6!=2^{10}.3^2.5$$
$$|W_6^\prime|=|W_6^+|=|\WC_6^\prime|=23040=2^5.6!=2^9.3^2.5.\leqno{\mathrm{et}}$$

\bigskip

\lemma{irreductible}
{\it Le groupe $D(W_6)$ est un sous-groupe irr\'eductible de 
$\Gb\Lb(\www V)$.}

\bigskip

\proof Supposons trouv\'e un sous-espace vectoriel $V'$ non trivial de $\www V$ qui 
est stable par $D(W_6)$. Par semi-simplicit\'e, on peut supposer que $\dim V' \le 3$. 
Puisque $V'$ est stable sous l'action de $\AG_6$, on a 
$V'=\CM(\e_1+\e_2+\e_3+\e_4+\e_5+\e_6)$. Mais 
$\diag(-1,-1,1,1,1,1) \in A_6^\prime$ ne stabilise alors pas $V'$, ce 
qui est contraire \`a l'hypoth\`ese.\fin

\bigskip

\lemma{engendre}
{\it Soit $H$ un sous-groupe distingu\'e de $W_6^\prime$ (\resp $W_6^+$, 
\resp $\WC_6^\prime$) tel que $\pi(H)=\SG_6$. Alors $H$ est \'egal \`a 
$W_6^\prime$ (\resp $W_6^+$, \resp $\WC_6^\prime$).}

\bigskip

\proof Il nous suffit de montrer que $A_6^\prime \incl H$. 
Si $\s \in \SG_6$, on note $\sigt$ un \'el\'ement de $H$ tel 
que $\pi(\sigt)=\s$. Si $a \in A_6^\prime$, alors 
$\sigt a \sigt^{-1} =\s a \s^{-1}$. D'autre part, 
$\sigt a \sigt^{-1} a^{-1} = [\s,a] \in H$ car $H$ est 
distingu\'e. Le r\'esultat d\'ecoule alors de \ref{eng}.\fin

\bigskip

\soussection{Une d\'efinition de $G_{31}$\label{unc}} 
Notons $\ssss(V)$ le sous-groupe de $\Gb\Lb(V)$ 
form\'e des \'el\'ements de d\'eterminant $1$ ou $-1$ et posons 
$\ISO(V,\wgras)=<\Sb\Ob(V,\wgras),i \Id_{\www V} >$. 
Il est alors clair que $\ssss(V) = < \Sb\Lb(V),\z \Id_V >$ 
(o\`u $\z \in \CM^\times$ est une racine primitive huiti\`eme de l'unit\'e). 
Par cons\'equent, 
$$\L^{-1}(\ISO(\www V,\wgras))=\ssss(V)\qquad{\mathrm{et}}\qquad
\L(\ssss(V))= \ISO(\www V,\wgras).$$
Nous posons la d\'efinition suivante~:

\bigskip

\begin{centerline}
{\begin{tabular}{|c|}
\hline
${\vphantom{{\DS{A_{\DS{B}}^{\DS{A}}}} \over {\DS{A_{\DS{B}}^{\DS{A}}}}}} 
\quad G_{31}=\L^{-1}(\WC_6^\prime)\quad $ \\
\hline
\end{tabular}}
\end{centerline}

\bigskip

\noindent 
Alors $G_{31}$ est un sous-groupe de $\ssss(V)$ car $\WC_6^\prime \incl \ISO(V,\wgras)$. 
Avant de construire un ensemble de r\'eflexions engendrant $G_{31}$, nous 
\'enon\c{c}ons un r\'esultat qui nous sera utile~:

\bigskip

\lemma{pm 1}
{\it Soit $G$ un sous-groupe de $G_{31}$ tel que $\L(G)=\WC_6^\prime$. 
Alors $G=G_{31}$.}

\bigskip

\proof On a $-\Id_{\www V} \in \WC_6^\prime$. Donc au moins l'un 
des \'el\'ements $i \Id_V$ ou $-i\Id_V$ appartient \`a $G$. 
Donc $-\Id_V=(i\Id_V)^2=(-i\Id_V)^2 \in G$. Par cons\'equent, $G=G_{31}$.\fin

\bigskip

Notons $\m_0$ l'\'el\'ement de $W_6^+$ dont la matrice dans 
la base $\BC$ est
$$\matrice{0 & 1 & 0 & 0 & 0 & 0 \\
          -1 & 0 & 0 & 0 & 0 & 0 \\
	   0 & 0 & 0 & 1 & 0 & 0 \\
	   0 & 0 &-1 & 0 & 0 & 0 \\
	   0 & 0 & 0 & 0 & 0 & 1 \\
	   0 & 0 & 0 & 0 &-1 & 0}.$$
On note $M$ la classe de conjugaison de $\m_0$ dans $W_6^+$. Notons que 
$\m_0$ et $-\m_0$ sont conjugu\'es dans $W_6$ mais ne le sont pas 
dans $W_6^+$~: en effet, on a 
$$-\m_0=t_0 \m_0 t_0^{-1},$$
o\`u $t_0=\diag(1,-1,1,-1,1,-1) \not\in W_6^+$ et un calcul \'el\'ementaire 
donne que $C_{W_6}(\m_0) \incl W_6^+$. Cela montre que la classe de 
conjugaison de $\m_0$ dans $W_6$ est la r\'eunion disjointe de 
$M$ et $-M$. 

Si $\m \in M$, alors le polyn\^ome caract\'eristique de $\m$ est 
$(X-i)^3(X+i)^3$. Par suite, le polyn\^ome caract\'eristique de $i\m$ est 
$(X-1)^3(X+1)^3$. 
Donc $\L^{-1}(\{i\m,-i\m\})=\{\mut,i\mut,-\mut,-i\mut\}$ o\`u $\mut$ est 
une r\'eflexion. Or, $i\m$ et $-i\m$ appartiennent \`a $\WC_6^\prime$ 
donc $\mut \in G_{31}$. On pose $\Mti=\{\mut~|~\m \in M\}$~; 
$\Mti$ est un ensemble de r\'eflexions de $G_{31}$. 

\medskip

{\it Par la suite, 
on supposera que $i$ est choisi de sorte que $\L^{-1}(i\m_0)=\{\mut_0,-\mut_0\}$.} 

\bigskip

\lemma{conju mu}
{\it Si $\m \in M$, alors $\L^{-1}(i\m)=\{\mut,-\mut\}$. 
De plus, $\Mti$ est une classe de conjugaison de $G_{31}$.}

\bigskip

\proof  Soit $\m \in M$. Puisque $\m$ et $\m_0$ sont conjugu\'es dans $W_6^+$ 
et puisque $\m_0 \not\in D(W_6)$, il existe $w \in D(W_6)$ tel que 
$\m=w\m_0w^{-1}$. De plus, $D(W_6) \incl \WC_6^\prime$, donc il existe 
$g \in G_{31}$ tel que $\L(g)=w$. Par suite, 
$\L(g \mut_0 g^{-1}) = w i\m_0 w^{-1}=i\m$. Puisque $g \mut_0 g^{-1}$ 
est une r\'eflexion, on en d\'eduit que $g \mut_0 g^{-1} = \mut$. 
Donc $\L^{-1}(i\m)=\{\mut,-\mut\}$ et $\Mti$ est bien une classe 
de conjugaison dans $G_{31}$.\finl

\bigskip

\theorem{main}
{\it Le groupe $G_{31}$ est un sous-groupe de r\'eflexion complexe 
irr\'eductible d'ordre $46080$ de $\Gb\Lb(V)$. On a $G_{31}=<\Mti>$.}

\bigskip

\proof Tout d'abord $G_{31}$ est irr\'eductible car son image $\WC_6^\prime$ 
est un sous-groupe \irr de $\Gb\Lb(\www V)$ d'apr\`es le lemme \ref{irreductible}. 
D'autre part, $|G_{31}|=2|\WC_6^\prime|=|W_6|$. 

Notons $G$ le sous-groupe de $G_{31}$ engendr\'e par $\Mti$. 
Alors $\L(G)$ est le sous-groupe de $\WC_6^\prime$ engendr\'e par $\L(\Mti)$. 
Or, d'apr\`es le lemme \ref{conju mu}, $\L(\Mti)$ 
est une classe de conjugaison de $\WC_6^\prime$. 
Donc $G$ un sous-groupe distingu\'e de $\WC_6^\prime$. 
D'autre part, $\pi(\L(G))$ est un sous-groupe distingu\'e de $\SG_6$ 
contenant $\pi(\m_0)=(1,2)(3,4)(5,6)$. Donc $\pi(\L(G))=\SG_6$. 
Il r\'esulte du lemme \ref{engendre} que $\L(G)=\WC_6^\prime$. 
Donc, d'apr\`es le lemme \ref{pm 1}, on a $G=G_{31}$.\fin

\bigskip

Le th\'eor\`eme \ref{main} montre que le 
groupe $G_{31}$ est bien le groupe de r\'eflexion complexe 
not\'e $G_{31}$ dans la classification de Shephard et Todd \cite{ST}. 

\bigskip

\section{Quelques propri\'et\'es du groupe $G_{31}$}~

\medskip

Nous allons ici \'etablir quelques propri\'et\'es bien connues du groupe $G_{31}$ 
(groupe d\'eriv\'e, centre, $2$-sous-groupes distingu\'es, nombre minimal 
de r\'eflexions engendrant $G_{31}$) en utilisant simplement 
les propri\'et\'es du groupe $\WC_6^\prime$. 
Dans cette section, nous noterons $\r$ la matrice 
$$\matrice{-1 & 0 & 0 & 0 & 0 & 0 \\ 0 & 1 & 0 & 0 & 0 & 0 \\ 0 & 0 & 1 & 0 & 0 & 0 \\
            0 & 0 & 0 & 1 & 0 & 0 \\ 0 & 0 & 0 & 0 & 0 & 1 \\ 0 & 0 & 0 & 0 & 1 & 0}.$$
Alors $\r \in W_6^+$ et $\g(\r)=-1$. Donc 
\equat\label{irrationnel}
i\r \in \WC_6^\prime \qquad{\mathrm{et}}\qquad \Tr(i\r)=2i \not\in \QM.
\endequat

\bigskip

\soussection{Scission des suites exactes} Nous allons 
ici essentiellement \'etudier les suites exactes \ref{w} \`a \ref{dw}. 
Nous allons notamment d\'eterminer si elles sont scind\'ees ou non. 
Nous aurons besoin pour cela du r\'esultat \'el\'ementaire suivant~:

\bigskip

\lemma{t}
{\it Soit $f : \SG_6 \to \Tb \rtimes \SG_6$ un morphisme de groupes tel que 
$\pi \circ f = \Id_{\SG_6}$. Alors $f(\SG_6) \incl \Ker \g$.}

\bigskip

\proof Notons $c=(1,2,3,4,5,6) \in \SG_6$. Si 
$(t_1,t_2,t_3,t_4,t_5,t_6) \in (\CM^\times)^6$, alors
$$c \diag(t_1,t_2,t_3,t_4,t_5,t_6) c^{-1}= \diag(t_6,t_1,t_2,t_3,t_4,t_5).$$
\'Ecrivons 
$$f(c)=\diag(a_1,a_2,a_3,a_4,a_5,a_6) c$$
avec $(a_1,a_2,a_3,a_4,a_5,a_6) \in (\CM^\times)^6$. Puisque $f(c)^6=\Id_{\www V}$, 
on a $a_1a_2a_3a_4a_5a_6 = 1$. Donc $f(c) \in \Ker \g$. Par cons\'equent, 
$f(c') \in \Ker \g$ pour tout cycle $c'$ de longueur $6$. Donc 
$f(\SG_6) \incl \Ker \g$.\fin

\bigskip

Pour la suite de cet article, nous n'aurons besoin que des r\'esultats 
de la proposition suivante concernant les groupes $W_6^\prime$ et $\WC_6^\prime$. 
Nous \'enon\c{c}ons cependant les autres par souci d'exhaustivit\'e.

\bigskip

\proposition{scindage}
{\it $({\mathrm{a}})$ 
Les suites \ref{w}, \ref{w'} et \ref{dw} sont scind\'ees gr\^ace 
aux inclusions 
$\SG_6 \incl W_6^\prime \incl W_6$ et $\AG_6 \incl D(W_6)$.

\tete{b} La suite exacte 
$$1 \longto A_6^\prime/Z(\WC_6^\prime) \longto \WC_6^\prime/Z(\WC_6^\prime) 
\longto \SG_6 \longto 1$$ 
n'est pas scind\'ee. 
Les groupes $\WC_6^\prime/Z(\WC_6^\prime)$ et $W_6^\prime/Z(W_6^\prime)$ 
ne sont pas isomorphes.

\tete{c} Les suites exactes \ref{w+}, \ref{wc} et \ref{wc'} ne sont pas scind\'ees.

\tete{d} Les groupes $W_6$ et $\WC_6$ ne sont pas isomorphes. 

\tete{e} Les groupes 
$\WC_6^\prime$ et $W_6^\prime$ ne sont pas isomorphes.

\tete{f} Les groupes $W_6^+$ et $W_6^\prime$
ne sont pas isomorphes.}

\bigskip

\noindent{\sc Remarque - } Bien qu'ils ne soient pas isomorphes, un calcul 
avec {\tt GAP} montre que les groupes $W_6^\prime/Z(W_6^\prime)$ et 
$\WC_6^\prime/Z(\WC_6^\prime)$ ont la m\^eme table de caract\`eres. 
En revanche, les groupes $W_6^\prime$ et $\WC_6^\prime$ ont tous deux 
$37$ classes de conjugaison mais n'ont pas la m\^eme table de caract\`eres~: 
celle de $\WC_6^\prime$ contient des valeurs irrationnelles comme 
le montre \ref{irrationnel}.

Pour finir, toujours d'apr\`es {\tt GAP}, les groupes $W_6^+$ et 
$\WC_6^\prime$ ne sont pas isomorphes.\finl

\bigskip

\proof (a) est \'evident. 

\medskip

\tete{b} Compte tenu du (a) et du fait que 
$A_6^\prime=O_2(\WC_6^\prime)=O_2(W_6^\prime)$, les deux assertions sont 
\'equivalentes. Nous allons montrer que la suite exacte n'est pas scind\'ee. 

Supposons-l\`a scind\'ee. Alors il existe deux \'el\'ements $g$ 
et $h$ dans $\WC_6^\prime$ tels que 

\begin{quotation}
{$(\a)$ \qquad $\pi(g)=(1,2,3,4)$ et $\pi(h)=(5,6)$.

\noindent$(\b)$ \qquad $h^2 \in Z(\WC_6^\prime)$.

\noindent$(\g)$ \qquad $[g,h] \in Z(\WC_6^\prime)$.}
\end{quotation}

\noindent D'apr\`es \ref{irrationnel} et $(\a)$, il existe $a \in A_6^\prime$ tel que 
$h=ia\r$. Donc, d'apr\`es $(\b)$, on a $a=\diag(a_1,a_2,a_3,a_4,b,b)$ 
avec $a_1^2=a_2^2=a_3^2=a_4^2=b^2=a_1a_2a_3a_4=1$. D'autre part, d'apr\`es $(\a)$, 
il existe $a' \in A_6^\prime$ tel que $g=ia'\s$, o\`u 
$$\s=\matrice{0 & 0 & 0 & 1 & 0 & 0 \\
              1 & 0 & 0 & 0 & 0 & 0 \\
	      0 & 1 & 0 & 0 & 0 & 0 \\
	      0 & 0 & 1 & 0 & 0 & 0 \\
	      0 & 0 & 0 & 0 & 1 & 0 \\
	      0 & 0 & 0 & 0 & 0 &-1 \\}.$$
Si on pose $a'=\diag(a_1^\prime,a_2^\prime,a_3^\prime,a_4^\prime,a_5^\prime,a_6^\prime)$, 
alors 
$$[g,h]=
\diag(-a_1a_4,-a_1a_2,a_2a_3,a_3a_4,-a_5^\prime a_6^\prime,-a_5^\prime a_6^\prime).$$
Donc, d'apr\`es $(\g)$, on a $-a_1a_4=a_2a_3$, et donc $a_1a_2a_3a_4=-1$, 
ce qui est impossible.

\medskip

\tete{c} D'apr\`es (b), la suite \ref{wc'} n'est pas scind\'ee. 
Puisque $W_6^+ \incl \WC_6$, il ne reste plus qu'\`a montrer que 
la suite \ref{wc} n'est pas scind\'ee. 
Supposons qu'elle est scind\'ee. Notons $f : \SG_6 \to \WC_6$ 
une section du morphisme $\pi$. D'apr\`es le lemme \ref{t}, $f(\SG_6) \incl \Ker \g$. 
Donc $f(\SG_6) \incl \WC_6^\prime$. En d'autres termes, la suite \ref{wc'} 
est scind\'ee, ce qui n'est pas possible. 

\medskip

\tete{d} D'apr\`es le lemme \ref{irreductible}, 
le centre de $\WC_6$ est d'ordre $4$ tandis que celui de $W_6$ est 
d'ordre $2$.

\medskip

\tete{e} et (f) d\'ecoulent imm\'ediatement de (a), (c) et (d).\fin

\bigskip

\soussection{Structure du groupe $G_{31}$} 
La proposition suivante se d\'eduit des propri\'et\'es du groupe 
$\WC_6^\prime$ d\'emontr\'ees pr\'ec\'edemment. 

\bigskip

\proposition{centre derive}
{\it $({\mathrm{a}})$ $\Mti$ est l'unique classe de conjugaison de $G_{31}$ 
form\'ee de r\'eflexions. 

\tete{b} $D(G_{31})=G_{31} \cap \Sb\Lb(V)$ et $|G_{31}/D(G_{31})|=2$.

\tete{c} $Z(G_{31}) = <i \Id_V >$ et $Z(G_{31}) \incl D(G_{31})$.

\tete{d} $\L$ induit un isomorphisme $D(G_{31})/Z \simeq D(W_6)$.

\tete{e} $O_2(G_{31})=\L^{-1}(A_6^\prime)\incl D(G_{31})$ est non ab\'elien et 
$O_2(G_{31})/Z \simeq A_6^\prime \simeq (\ZM/2\ZM)^5$.

\tete{f} Les $2$-sous-groupes distingu\'es de $G_{31}$ sont $\{1\}$, 
$Z$, $Z(G_{31})$ et $O_2(G_{31})$.

\tete{g} Les suites exactes 
$$\diagram
1 \rrto &&O_2(G_{31})/Z(G_{31}) \rrto &&G_{31}/Z(G_{31}) 
\rrto^{\DS{\pi \circ \L}} &&\SG_6 \rrto &&1
\enddiagram$$
$$\diagram
1 \rrto &&O_2(G_{31}) \rrto &&G_{31} \rrto^{\DS{\pi \circ \L}} &&\SG_6 
\rrto &&1
\enddiagram\leqno{\mathit{et}}$$
ne sont pas scind\'ees.

\tete{h} Les groupes $G_{31}/Z(G_{31})$ et $W_6^\prime/Z(W_6^\prime)$ 
ne sont pas isomorphes.}

\bigskip

\noindent{\sc Remarque - } Dans \cite[table 3]{BMR}, il est annonc\'e 
que les groupes $G_{31}/Z(G_{31})$ et $W_6^\prime/Z(W_6^\prime)$ 
ne sont pas isomorphes (ce qui est exactement le (h) de la proposition  
\ref{centre derive} ci-dessus). En revanche, toujours dans \cite[table 3]{BMR}, 
il est dit que la premi\`ere suite exacte du (g) est scind\'ee~: ce n'est pas 
vrai comme le prouve l'argument ci-dessous (et comme 
peuvent le montrer aussi des calculs effectu\'es avec {\tt GAP}).\finl

\bigskip

\proof (a) Soit $s \in G_{31}$ une r\'eflexion. Alors $-i\L(s) \in W_6^+$ 
a pour polyn\^ome caract\'eristique $(X-i)^3(X+i)^3$. 
Or, les seuls \'el\'ements de $W_6^+$ dont le polyn\^ome caract\'eristique 
est $(X-i)^3(X+i)^3$ sont les conjugu\'es sous $W_6$ 
de $\m_0$. Donc $-i\L(s)=\e\m$ avec $\e \in \{1,-1\}$ et $\m \in M$. 
Donc $s=\mut$ car $s$ est une r\'eflexion. 

\medskip

\tete{b} D'apr\`es (a), $G_{31}$ est engendr\'e par une classe de conjugaison 
d'\'el\'ements d'ordre $2$ donc $|G_{31}/D(G_{31})| \le 2$. Mais 
$D(G_{31})$ est contenu dans $G_{31} \cap \Sb\Lb(V)$ et ce dernier 
est d'indice $2$ dans $G_{31}$. D'o\`u le r\'esultat.

\medskip

\tete{c} d\'ecoule de (b) et du fait que $G_{31}$ est irr\'eductible sur $V$ 
(voir theor\`eme \ref{main}).

\medskip

\tete{d} Il suffit de remarquer que $Z \incl D(G_{31})$, ce qui r\'esulte 
de (b).

\medskip

\tete{e} Le fait que $O_2(G_{31})=\L^{-1}(A_6^\prime)$ et que 
$O_2(G_{31})/Z \simeq A_6^\prime$ r\'esulte de ce que $O_2(\WC_6^\prime)=A_6^\prime$ 
car $O_2(\SG_6)=1$. Le fait que $O_2(G_{31}) \incl D(G_{31})$ d\'ecoule 
alors du (b). 
Il nous reste \`a montrer que $O_2(G_{31})$ 
n'est pas ab\'elien. 

Le nombre minimal de g\'en\'erateurs de $O_2(G_{31})$ est $5$ car 
$O_2(G_{31})/Z \simeq A_6^\prime \simeq (\ZM/2\ZM)^5$. Or, un groupe ab\'elien dont 
le nombre minimal de
g\'en\'erateurs est $5$ n'admet pas de repr\'esentation fid\`ele de 
dimension $4$. Donc $O_2(G_{31})$ n'est pas ab\'elien. 

\medskip

\tete{f} Soit $H$ un $2$-sous-groupe distingu\'e de $G_{31}$ diff\'erent 
de $\{1\}$, $Z$ et $O_2(G_{31})$. Alors $Z \incl H$ et, en appliquant 
le morphisme $\L$, on obtient que $\L(H)$ est un $2$-sous-groupe 
non trivial de $A_6^\prime$ stable par l'action de $\SG_6$. 
Donc $\L(H)=\{\Id_{\www V},-\Id_{\www V} \}$, ce qui montre que 
$H=Z(G_{31})$.

\medskip

\tete{g} Il est clair que la premi\`ere suite est exacte. Le fait 
qu'elle est non scind\'ee d\'ecoule de la proposition \ref{scindage} (b). 
Il en r\'esulte imm\'ediatement que la deuxi\`eme suite exacte 
n'est pas scind\'ee. 

\medskip

\tete{h} Puisque $G_{31}/Z(G_{31})=\WC_6^\prime/Z(\WC_6^\prime)$, le r\'esultat 
d\'ecoule de la proposition \ref{scindage} (b).\fin

\bigskip

\soussection{Engendrement par $5$ r\'eflexions} 
La proposition suivante pr\'ecise le lemme \ref{pm 1}. 

\bigskip

\proposition{max}
{\it Si $G$ est un sous-groupe de $G_{31}$, alors les assertions suivantes 
sont \'equivalentes~:

\tete{1} $G_{31}=G$.

\tete{2} $\WC_6^\prime=\L(G)$.

\tete{3} $\SG_6=\pi(\L(G))$.}

\bigskip

\proof L'\'equivalence entre (1) et (2) d\'ecoule du lemme \ref{pm 1}. 
Il est clair que (1) implique (3). Montrons que (3) implique (1). 
On suppose donc que $\SG_6=\pi(\L(G))$. On pose
$A=O_2(G_{31}) \cap G$ et $A'=AZ$. 

Puisque $Z \incl A'$ et $O_2(G_{31})/Z$ est ab\'elien, $A'$ est distingu\'e 
dans $O_2(G_{31})$. D'autre part, il est normalis\'e par $G$. Donc 
il est distingu\'e dans $G_{31}$. Donc, d'apr\`es la proposition 
\ref{centre derive} (f), on a $A'=Z$, $Z(G_{31})$ ou $O_2(G_{31})$. 
Si $A'=Z$ ou $Z(G_{31})$, alors 
la premi\`ere suite exacte de la proposition \ref{centre derive} (g) 
serait scind\'ee, ce qui est impossible. 
Donc $A'=O_2(G_{31})$. Par cons\'equent, $G.Z=GA'=G_{31}$ et donc $G=G_{31}$ 
d'apr\`es le lemme \ref{pm 1}.\fin

\bigskip

\proposition{5}
{\it Le groupe $G_{31}$ est engendr\'e par $5$ r\'eflexions mais ne peut 
pas \^etre engendr\'e par $4$ r\'eflexions.}

\bigskip

\proof Nous noterons 
$\t : \SG_6 \to \SG_6$ un automorphisme non int\'erieur de $\SG_6$. 
Par construction, $\pi(M)=\pi(\L(\Mti))$ 
est la classe de conjugaison de $\SG_6$ form\'ee 
des produits de trois transpositions \`a supports disjoints. Par cons\'equent, 
$\t(\pi(M))$ est la classe de conjugaison de $\SG_6$ form\'e des 
transpositions (\cad des r\'eflexions). 

Si $1 \le j \le 5$, on pose $s_j=(j,j+1)$ et $w_j=(1,j)(2,j+1)$, ce qui  
implique que $s_j=w_j s_1 w_j^{-1}$. Posons maintenant 
$\m_j=\t^{-1}(w_j)\m_0 \t^{-1}(w_j)^{-1}$. 
Alors $\pi(\m_j)=\t^{-1}(s_j)$ et donc    
$\SG_6=<\pi(\m_1), \pi(\m_2),\pi(\m_3), \pi(\m_4), \pi(\m_5) >$. 
Par cons\'equent, d'apr\`es la proposition \ref{max}, on a 
$G_{31}=<\mut_1,\mut_2,\mut_3,\mut_4,\mut_5>$. Cela montre la premi\`ere assertion. 

Montrons la deuxi\`eme. Soit $\Sti$ un ensemble de r\'eflexions 
qui engendrent $G_{31}$. Alors $\pi(\L(\Sti))$ engendre 
$\SG_6$ et donc $\t(\pi(\L(\Sti)))$ engendre $\SG_6$. 
Or il est bien connu que $\SG_6$ ne peut pas \^etre 
engendr\'e par $4$ transpositions (en effet, un sous-groupe de 
$\SG_6$ engendr\'e par $4$ transpositions ne peut pas \^etre 
transitif sur l'ensemble $\{1,2,3,4,5,6\}$). Donc $|\Sti| \ge 5$.\fin

\end{document}